\documentclass{amsart}

\usepackage{fancyhdr}
\usepackage{lastpage}
\usepackage{stmaryrd,yhmath}

\newcommand{\bdis}{\begin{displaymath}}
\newcommand{\edis}{\end{displaymath}}
\newcommand{\be}{\begin{equation}}
\newcommand{\ee}{\end{equation}}

\newcommand{\mcal}{\mathcal}

\newcommand{\vt}{\vartheta}

\newcommand{\tg}{\tilde{\gamma}}

\newcommand{\zf}{\zeta\left(\frac{1}{2}+it\right)}

\newtheorem{lemma}[]{Lemma}

\theoremstyle{definition}

\theoremstyle{remark}

\newtheorem*{mydef1}{{\bf Theorem}}

\newtheorem*{mydef8}{{\bf Remark}}

\numberwithin{equation}{section}

\begin{document}

\title{Distribution of the roots of the equations $Z(t)=0$, $Z'(t)=0$ in the theory of the Riemann zeta-function}

\author{Jan Moser}

\address{Department of Mathematical Analysis and Numerical Mathematics, Comenius University, Mlynska Dolina M105, 842 48 Bratislava, SLOVAKIA}

\email{jan.mozer@fmph.uniba.sk}


\begin{abstract}
Let the symbols $\{\gamma\},\ \{t_0\};\ t_0\not=\gamma$ denote the sequences of the roots of the equations
$$Z(t)=0,\quad \text{and}\qquad  Z'(t)=0,$$
respectively, and
$$m(t_0)=\min\{ \gamma''-t_0,t_0-\gamma'\},\quad Q(t_0)=\max\{\gamma''-t_0,t_0-\gamma'\},\quad \gamma'<t_0<\gamma'',$$
where $\gamma',\gamma''$ are the neighboring zeroes. We have proved the following in this paper: on the Riemann hypothesis we have
$$\frac{Q(t_0)}{m(t_0)}<t_0\ln^2t_0\ln_2t_0\ln_3t_0,\quad t_0\to\infty.$$
This paper is the English version of the paper of ref. \cite{5}.
\end{abstract}

\maketitle

\section{Formulation of Theorem}

We will denote by the symbols
\bdis
\{\gamma\},\quad \{ t_0\}
\edis
the sequences of the roots of the equations
\bdis
Z(t)=0,\quad Z'(t)=0,\quad t_0\not=\gamma,
\edis
where (see \cite{6}, pp. 79, 329)
\be \label{1.1}
\begin{split}
& Z(t)=e^{i\vt(t)}\zf, \\
& \vt(t)=-\frac t2\ln\pi+\text{Im}\ln\Gamma\left(\frac 14+\frac 12it\right)= \\
& =\frac t2\ln\frac{t}{2\pi}-\frac t2-\frac{\pi}{8}+\mcal{O}\left(\frac 1t\right).
\end{split}
\ee

\begin{mydef8}
On the Riemann hypothesis the points of the sequences $\{\gamma\},\{ t_0\}$ are separated each from other, i. e.
\bdis
\gamma'< t_0 <\gamma'',
\edis
where $\gamma',\gamma''$ are neighboring points of the sequence $\{\gamma\}$ (see \cite{4}, Corollary 3).
\end{mydef8}

The sequence $\{ t_0\}$ itself oscillates in a complicated manner around the sequence $\{\gamma\}$, (comp., for example, the graph of the
function $Z(t)$ in the neighborhood of the first Lehmer pair of zeroes in \cite{1}, p. 296). Let further
\be \label{1.2}
\begin{split}
& Q(t_0)=\max\{\gamma''-t_0,t_0-\gamma'\}, \\
& m(t_0)=\min\{ \gamma''-t_0,t_0-\gamma'\}.
\end{split}
\ee
Then the quotient
\bdis
\frac{Q(t_0)}{m(t_0)}
\edis
measures the asymmetry of the point $t_0$ with respect to the points $\gamma',\gamma''$. \\

In this paper we prove the following statement.

\begin{mydef1}
On the Riemann hypothesis we have
\be \label{1.3}
\frac{Q(t_0)}{m(t_0)}<t_0\ln^2t_0\ln_2t_0\ln_3t_0,\quad t_0\to\infty,
\ee
where
\bdis
\ln_2t_0=\ln\ln t_0,\dots
\edis
\end{mydef1}

The estimate (\ref{1.3}) follows from our formula
\be \label{1.4}
\frac{\pi}{4}=\sum_{\gamma}\frac{t_0}{\gamma^2-t_0^2}+\mcal{O}\left(\frac{1}{t_0}\right)
\ee
that is a consequence of the Riemann hypothesis (see \cite{3}, p. 117, comp. \cite{4}). The formula (\ref{1.4}) is a conjugate formula
with respect to the Riemann formula
\bdis
c+2-\ln 4\pi=\sum_{\gamma}\frac{1}{\frac 14+\gamma^2},
\edis
where $c$ stands for the Euler's constant. Let us remind, finally, that on the Riemann hypothesis also the Littlewood estimate
\be \label{1.5}
\gamma''-\gamma'<\frac{A}{\ln\ln\gamma'},\quad \gamma'\to\infty
\ee
holds true.

\section{Lemmas}

\begin{lemma}
\be \label{2.1}
\sum_{t_0+T\leq\gamma}\frac{1}{\gamma^2-t_0^2}<A\frac{\ln T}{T},
\ee
where
\be \label{2.2}
T=T(t_0)=t_0\ln t_0\ln_2t_0.
\ee
\end{lemma}

\begin{proof}
Let $N(t)$ denote the number of zeroes
$$\beta+i\gamma$$
of the function
$$\zeta(s),\quad s=\sigma+it,\quad \beta\in (0,1),\ \gamma\in (0,t].$$
Since (see \cite{6}, p. 178)
\bdis
N(t+1)-N(t)<A\ln t,
\edis
then
\bdis
\begin{split}
& \sum_{t_0+T\leq\gamma}\frac{1}{\gamma^2-t_0^2}=\sum_{n=1}^\infty \sum_{\{t_0+T+n-1<\gamma\leq t_0+T+n\}}\frac{1}{\gamma^2-t_0^2}< \\
& < A\sum_{n=1}^\infty\frac{\ln(n+T-1+t_0)}{(n+T-1)(n+T-1-2t_0)}< \\
& < A\frac{\ln(T+t_0)}{T(T+2t_0)}+A\int_1^\infty\frac{\ln (x+T-1+t_0)}{(x+T-1)(x+T-1+2t_0)}{\rm d}x.
\end{split}
\edis
Following the identity
\bdis
\int_a^\infty \tau e^{-\tau}{\rm d}\tau=e^{-a}(a+1)
\edis
we obtain
\bdis
\begin{split}
& \int_1^\infty\frac{\ln (x+T-1+t_0)}{(x+T-1)(x+T-1+2t_0)}{\rm d}x<\int_1^\infty\frac{\ln (x+T-1+t_0)}{(x+T-1)(x+T-1+t_0)}{\rm d}x= \\
& =\int_{\ln(T+t_0)}^\infty\frac{\tau}{e^\tau-t_0}{\rm d}\tau=\int_{\ln(T+t_0)}^\infty\frac{e^{-\tau}\tau}{1-t_0e^{-\tau}}{\rm d}\tau<\frac{T+t_0}{T}\int_{\ln(T+t_0)}^\infty e^{-\tau}\tau{\rm d}\tau= \\
& = \frac 1T\left\{ \ln(T+t_0)+1\right\}<A\frac{\ln T}{T},
\end{split}
\edis
i. e. we have finished the proof.
\end{proof}

Consequently, we have from (\ref{1.4}) by (\ref{2.1}), (\ref{2.2}) the following

\begin{lemma}
On the Riemann hypothesis we have
\be \label{2.3}
\begin{split}
& \frac{\pi}{8t_0}=\sum_{0<\gamma<t_0+T}\frac{1}{\gamma^2-t_0^2}+\mcal{O}\left(\frac{1}{t_0\ln_2t_0}\right),\quad t_0\to\infty, \\
& T=T(t_0)=t_0\ln t_0\ln_2t_0.
\end{split}
\ee
\end{lemma}

\section{Proof of the Theorem}

If there is any $\tilde{t}_0$ such that
\be \label{3.1}
K(\tilde{t}_0)=\frac{Q(\tilde{t}_0)}{m(\tilde{t}_0)}\geq \tilde{t}_0\ln^2\tilde{t}_0\ln_2\tilde{t}_0\ln_3\tilde{t}_0
\ee
then we have:

\subsection*{(A)}

If
\be \label{3.2}
\tilde{t}_0-\tilde{\gamma}'=K(\tilde{t}_0)(\tilde{\gamma}''-\tilde{t}_0)
\ee
then by (2.3) there is some $A>0$ such that
\be \label{3.3}
\frac{A}{\tilde{t}_0}> \sum_{0<\gamma\leq \tilde{t}_0+\tilde{T}}\frac{1}{\gamma^2-\tilde{t}_0^2},\quad \tilde{T}=T(\tilde{t}_0).
\ee
Now $\tilde{\gamma}'<\tilde{t}_0<\tilde{\gamma}''$ and the symbol $\qquad n(\gamma)$ denotes the order of the zero
$\frac 12+i\gamma$. We have
\be \label{3.4}
\begin{split}
& \sum_{0<\gamma\leq \tilde{t}_0+\tilde{T}}\frac{1}{\gamma^2-\tilde{t}_0^2}=\sum_{0<\gamma\leq\tilde{\gamma}'}\frac{1}{\gamma^2-\tilde{t}_0^2}+
\frac{n(\tilde{\gamma}'')}{(\tilde{\gamma}'')^2-\tilde{t}_0^2}+\sum_{\tilde{\gamma}''<\gamma\leq \tilde{t}_0+\tilde{T}}\frac{1}{\gamma^2-\tilde{t}_0^2} > \\
& > \sum_{0<\gamma<\tilde{\gamma}'}\frac{1}{\gamma^2-\tilde{t}_0^2}+\frac{1}{(\tilde{\gamma}'')^2-\tilde{t}_0^2}.
\end{split}
\ee
Since (see \cite{6}, p. 181)
\be \label{3.5}
N(t)<At\ln t,
\ee
then
\be \label{3.6}
\begin{split}
& \sum_{0<\gamma<\tilde{\gamma}'}\frac{1}{\gamma^2-\tilde{t}_0^2}=\sum_{0<\gamma<\tilde{t}_0}\frac{1}{\gamma^2-\tilde{t}_0^2}=
-\sum_{0<\gamma<\tilde{t}_0}\frac{1}{(\tilde{t}_0-\gamma)(\gamma+\tilde{t}_0)}> \\
& > -\frac{1}{\tilde{t}_0-\tilde{\gamma}'}\sum_{0<\gamma<\tilde{t}_0}\frac{1}{\gamma+\tilde{t}_0}>
-\frac{1}{\tilde{t}_0(\tilde{t}_0-\tilde{\gamma}')}\sum_{0<\gamma<\tilde{t}_0}1>\\
& > -\frac{1}{\tilde{t}_0(\tilde{t}_0-\tilde{\gamma}')}A\tilde{t}_0\ln \tilde{t}_0=-A\frac{\ln \tilde{t}_0}{\tilde{t}_0-\tilde{\gamma}'},
\end{split}
\ee
and by (3.2)
\be \label{3.7}
\frac{1}{(\tilde{\gamma}'')^2-\tilde{t}_0^2}=\frac{1}{(\tg''-\tilde{t}_0)(\tg''+\tilde{t}_0)}>\frac{1}{2\tg''(\tg''-\tilde{t}_0)}>
A\frac{K(\tilde{t}_0)}{\tilde{t}_0(\tilde{t}_0-\tg')}.
\ee
Consequently, from (3.4) by (1.5), (3.1), (3.6), (3.7) we have
\be \label{3.8}
\begin{split}
& \sum_{0<\gamma\leq \tilde{t}_0+\tilde{T}}\frac{1}{\gamma^2-\tilde{t}_0^2}>\frac{A}{\tilde{t}_0-\tg'}
\left\{\frac{K(\tilde{t}_0)}{\tilde{t}_0}-\ln\tilde{t}_0\right\}> \\
& > \frac{A}{\tg''-\tg'}\ln^2\tilde{t}_0\ln_2\tilde{t}_0\ln_3\tilde{t}_0>A\ln^2\tilde{t}_0\ln_2^2\tilde{t}_0\ln_3\tilde{t}_0,
\end{split}
\ee
and finally (see (3.3), (3.8))
\bdis
\frac{1}{\tilde{t}_0}>A\ln^2\tilde{t}_0\ln_2^2\tilde{t}_0\ln_3\tilde{t}_0,
\edis
i. e. we obtain the contradiction.

\subsection*{(B)}

If
\be \label{3.9}
\tg''-\tilde{t}_0=K(\tilde{t}_0)(\tilde{t}_0-\tg')
\ee
then by (2.3) there is some $A>0$ such that
\be \label{3.10}
\frac{A}{\tilde{t}_0}<\sum_{0<\gamma\leq \tilde{t}_0+\tilde{T}}\frac{1}{\gamma^2-\tilde{t}_0^2}.
\ee
In this case we have the following lower estimate for the corresponding sums: by (2.2) and (3.5) one obtains
\be \label{3.11}
\sum_{0<\gamma\leq \tilde{t}_0+\tilde{T}}1<A\tilde{T}\ln\tilde{T}<A\tilde{t}_0\ln^2\tilde{t}_0\ln_2\tilde{t}_0,
\ee
next, by (3.9)
\bdis
\sum_{0<\gamma\leq\tg'}\frac{1}{\gamma^2-\tilde{t}_0^2}<\frac{1}{(\tg')^2-\tilde{t}_0^2}=-\frac{K(\tilde{t}_0)}{(\tg'+\tilde{t}_0)(\tg''-\tilde{t}_0)}
<-A\frac{K(\tilde{t}_0)}{\tilde{t}_0(\tg'-\tilde{t}_0)},
\edis
and finally (see (3.11))
\bdis
\sum_{\tg''<\gamma\leq \tilde{t}_0+\tilde{T}}\frac{1}{\gamma^2-\tilde{t}_0^2}<\frac{1}{(\tg'')^2-\tilde{t}_0^2}\sum_{\tg''<\gamma\leq \tilde{t}_0+\tilde{T}}1
< A\frac{\ln^2\tilde{t}_0\ln_2\tilde{t}_0}{\tg''-\tilde{t}_0}.
\edis
Consequently,
\be \label{3.12}
\begin{split}
& \sum_{0<\gamma\leq\tilde{t}_0+\tilde{T}}\frac{1}{\gamma^2-\tilde{t}_0^2}=
\left\{\sum_{0<\gamma\leq\tg''}+\sum_{\tg''\leq\gamma\leq\tilde{t}_0+\tilde{T}}\right\}\frac{1}{\gamma^2-\tilde{t}_0^2}< \\
& < \frac{A}{\tg''-\tilde{t}_0}\left\{-\frac{K(\tilde{t}_0)}{\tilde{t}_0}+\ln^2\tilde{t}_0\ln_2\tilde{t}_0\right\}.
\end{split}
\ee
Hence, from (3.10) by (3.1), (3.12) we have that
\bdis
\frac{1}{\tilde{t}_0}<\frac{A}{\tg''-\tilde{t}_0}\{-\ln_3\tilde{t}_0+1\}\ln^2\tilde{t}_0\ln_2\tilde{t}_0<0,
\edis
i. e. we have obtained the contradiction. \\

\thanks{I would like to thank Michal Demetrian for helping me with the electronic version of this work.}

\end{document}